\documentclass{amsart}

\usepackage{amsmath} 
\usepackage{amssymb}

\newtheorem{theorem}{Theorem}[section] 
\newtheorem{claim}{Claim}[theorem]
\newtheorem{lemma}[theorem]{Lemma} 
\newtheorem{proposition}[theorem]{Proposition} 
\newtheorem{corollary}[theorem]{Corollary} 

\theoremstyle{definition}
\newtheorem{definition}[theorem]{Definition}
\newtheorem{example}[theorem]{Example}
\newtheorem{prob}{Problem}
\newtheorem{problem}[theorem]{Problem}

\theoremstyle{remark}
\newtheorem{remark}[theorem]{Remark}

\numberwithin{equation}{section}
\newcommand{\forces}{\Vdash}

\newcommand{\bV}{{\bf V}} 
\newcommand{\lesdot}{\mathrel{\mathord{<}\!\!\raise 
0.8 pt\hbox{$\scriptstyle\circ$}}}

\newcommand{\add}{{\rm {\bf add}}\/} 
 
\newcommand{\cof}{{\rm {\bf cof}}\/} 
\newcommand{\non}{{\rm {\bf non}}\/} 
\newcommand{\cov}{{\rm {\bf cov}}\/}
 
\newcommand{\gd}{{\mathfrak d\/}} 
\newcommand{\gb}{{\mathfrak b}}

\newcommand{\can}{2^{\textstyle \omega}} 
 
\newcommand{\baire}{\omega^{\textstyle \omega}} 
\newcommand{\iso}{[\omega]^{\textstyle \omega}}

\newcommand{\conc}{{}^\frown\!}
\newcommand{\lh}{{\rm lh}\/}
\newcommand{\rest}{{\restriction}}
\newcommand{\mrot}{{\rm root}\/} 
\newcommand{\suc}{{\rm succ}}

\newcommand{\rng}{{\rm Range}}

\newcommand{\nor}{{\rm {\bf nor}}\/} 
\newcommand{\pos}{{\rm pos}}

\newcommand{\dn}{{\rm dn}}
\newcommand{\up}{{\rm up}}

\newcommand{\FC}{{\rm FT}}
\newcommand{\IFC}{{\rm IFT}}

\newcommand{\tree}{{\rm tree}}
\newcommand{\rpH}{{\rm rp}_\bH}
\newcommand{\rpHn}{{\rm rp}_\bH^n}

\newcommand{\lev}{{\rm lev}}
\newcommand{\UM}{{\mathbb U}{\mathbb M}}
\newcommand{\bQp}{{\bbQ^\tree(\gp)}}

\newcommand{\cA}{{\mathcal A}}

\newcommand{\cB}{{\mathcal B}}

\newcommand{\bbC}{{\mathbb C}}

\newcommand{\bbD}{{\mathbb D}}

\newcommand{\bH}{{\bf H}}

\newcommand{\bF}{{\bf F}}
\newcommand{\cF}{{\mathcal F}}
\newcommand{\cG}{{\mathcal G}}
\newcommand{\cI}{{\mathcal I}}

\newcommand{\cM}{{\mathcal M}}

\newcommand{\bbP}{{\mathbb P}}

\newcommand{\gp}{{\mathfrak p}}
\newcommand{\bbQ}{{\mathbb Q}}
\newcommand{\dbQ}{{\name{\mathbb Q}}}
\newcommand{\mbR}{{\mathbb R}}

\newcommand{\cX}{{\mathcal X}}
\newcommand{\cY}{{\mathcal Y}}

\newcount\skewfactor
\def\mathunderaccent#1#2 {\let\theaccent#1\skewfactor#2
\mathpalette\putaccentunder}
\def\putaccentunder#1#2{\oalign{$#1#2$\crcr\hidewidth
\vbox to.2ex{\hbox{$#1\skew\skewfactor\theaccent{}$}\vss}\hidewidth}}
\def\name{\mathunderaccent\tilde-3 }

\begin{document}

\title{Universal forcing notions and ideals}

\author{Andrzej Ros{\l}anowski}
\address{Department of Mathematics\\
 University of Nebraska at Omaha\\
 Omaha, NE 68182-0243, USA}
\email{roslanow@member.ams.org}
\urladdr{http://www.unomaha.edu/$\sim$aroslano}
\thanks{The first author thanks the Hebrew University of Jerusalem for
 support during his visit to Jerusalem in Spring'2003. He also thanks his
 wife, Ma{\l}gorzata Jankowiak--Ros{\l}anowska for supporting him when he
 was working on this paper.}  

\author{Saharon Shelah}
\address{Einstein Institute of Mathematics\\
Edmond J. Safra Campus, Givat Ram\\
The Hebrew University of Jerusalem\\
Jerusalem, 91904, Israel\\
 and  Department of Mathematics\\
 Rutgers University\\
 New Brunswick, NJ 08854, USA}
\email{shelah@math.huji.ac.il}
\urladdr{http://www.math.rutgers.edu/$\sim$shelah}
\thanks{Both authors acknowledge support from the United States-Israel
  Binational Science Foundation (Grant no. 2002323). This is publication 845
  of the second author.}     

\subjclass{03E40, 03E17}
\date{April 2004}

\begin{abstract}
The main result of this paper is a partial answer to \cite[Problem
  5.5]{RoSh:672}: a finite iteration of Universal Meager forcing notions
adds generic filters for many forcing notions determined by universality
parameters. We also give some results concerning cardinal characteristics of
the $\sigma$--ideals determined by those universality parameters. 
\end{abstract}

\maketitle

\section{Introduction}
One of the most striking differences between measure and category was
discovered in Shelah \cite{Sh:176} where it was proved that the Lebesgue 
measurability of ${\boldsymbol\Sigma}^1_3$ sets implies $\omega_1$ is
inaccessible in ${\bf L}$, while one can construct (in ZFC) a forcing notion
$\bbP$ such that $\bV^\bbP\models$ ``projective subsets of $\mbR$ have the
Baire property''. For the latter result one builds a homogeneous ccc forcing
notion adding a lot of Cohen reals. Homogeneity is obtained by multiple use
of amalgamation (see \cite{JuRo} for a full explanation of how this works),
the Cohen reals come from compositions with the Universal Meager forcing
notion $\UM$ or with the Hechler forcing notion $\bbD$. The main point of
that construction was isolating a strong version of ccc, so called
{\em sweetness}, which is preserved in amalgamations.  Later, Stern
\cite{St85} introduced a weaker property, {\em topological sweetness}, which
is also preserved in amalgamations. Sweet (i.e., strong ccc) properties of
forcing notions were further investigated in \cite{RoSh:672}, where we
introduced a new property called {\em iterable sweetness} (see
\cite[Definition 4.2.1]{RoSh:672}) and we proved the following two results.  

\begin{theorem}
\begin{enumerate}
\item (See \cite[Proposition 4.2.2]{RoSh:672})
If $\bbP$ is a sweet ccc forcing notion (in the sense of \cite[Definition
  7.2]{Sh:176}) in which any two compatible elements have a least upper
bound, then $\bbP$ is iterably sweet.
\item (See \cite[Theorem 4.2.4]{RoSh:672})
If $\bbP$ is a topologically sweet forcing notion (in the sens of Stern
\cite[Definition 1.2]{St85}) and $\dbQ$ is a $\bbP$--name for an iterably
sweet forcing, then the composition $\bbP*\dbQ$ is topologically sweet. 
\end{enumerate}
\end{theorem}

In \cite[\S 2.3]{RoSh:672} we introduced a scheme of building forcing
notions from so called {\em universality parameters} (see \ref{univpar}
later). We proved that typically they are sweet (see \cite[Proposition
4.2.5]{RoSh:672}) and in natural cases also iterably sweet. So the question
arose if the use of those forcing notions in iterations gives us something
really new. Specifically, we asked:

\begin{prob}
[See {\cite[Problem 5.5]{RoSh:672}}]
\label{oldprob}
Is there a universality parameter $\gp$ satisfying the requirements of
\cite[Proposition 4.2.5(3)]{RoSh:672} such that no finite iteration of the
Universal Meager forcing notion adds a $\bQp$--generic real? Does the
Universal Meager forcing add generic reals for the forcing notions $\bQp$
defined from $\gp$ as in \ref{ucmz}, \ref{PPex}, \ref{PPpart} here? 
\end{prob}
 
Bad news is that Problem \ref{oldprob} has a partially negative answer: if the
universality parameter $\gp$ satisfies some mild conditions (i.e., is {\em
regular}, see \ref{regular}), then finite iteration of $\UM$ will add a
generic filter for the corresponding forcing notion, see Corollary
\ref{maincor}.   

Good news is that we have more examples of iterably sweet forcings, and 
they will be presented in a subsequent paper \cite{RoSh:856}.

The structure of the present paper is as follows. In the first section we
recall in a simplified form all the definitions and results we need from
\cite{RoSh:672}, and we define {\em regular universality parameters}.  
We also re-present the canonical examples we keep in mind in this
context. In the second section we prove our main result: a sequence Cohen
real --- dominating real --- Cohen real produces generic filters for forcing 
notions $\bQp$ determined by regular $\gp$ (see Theorem \ref{cohdomcoh}).   
In the following section we look at the $\sigma$--ideals $\cI_\gp$ for
regular $\gp$ and we prove a couple of inequalities concerning their
cardinal characteristics. 
\medskip

\noindent {\bf Notation}\quad Our notation is rather standard and compatible
with that of classical textbooks (like Jech \cite{J} or Bartoszy\'nski and
Judah \cite{BaJu95}). In forcing we keep the older convention that {\em a
  stronger condition is the larger one}. Our main conventions are listed
below.   

\begin{enumerate}
\item  For a forcing notion $\bbP$, all $\bbP$--names for objects in the
extension via $\bbP$ will be denoted with a tilde below (e.g.,
$\name{\tau}$, $\name{X}$). The complete Boolean algebra determined by
$\bbP$ is denoted by ${\rm RO}(\bbP)$.

\item For two sequences $\eta,\nu$ we write $\nu\vartriangleleft\eta$
whenever $\nu$ is a proper initial segment of $\eta$, and $\nu
\trianglelefteq\eta$ when either $\nu\vartriangleleft\eta$ or $\nu=\eta$. 
The length of a sequence $\eta$ is denoted by $\lh(\eta)$.

\item A {\em tree} is a family $T$ of finite sequences such that for some
$\mrot(T)\in T$ we have
\[(\forall\nu\in T)(\mrot(T)\trianglelefteq \nu)\quad\mbox{ and }\quad
\mrot(T)\trianglelefteq\nu\trianglelefteq\eta\in T\ \Rightarrow\ \nu\in T.\]
For a tree $T$, the family of all $\omega$--branches through $T$ is
denoted by $[T]$, and we let
\[\max(T)\stackrel{\rm def}{=}\{\nu\in T:\mbox{ there is no }\rho\in
T\mbox{ such that }\nu\vartriangleleft\rho\}.\]   
If $\eta$ is a node in the tree $T$ then 
\[\begin{array}{lcl}
\suc_T(\eta)&=&\{\nu\in T: \eta\vartriangleleft\nu\ \&\ \lh(\nu)=\lh(\eta)+1
\}\ \mbox{ and}\\
T^{[\eta]}&=&\{\nu\in T:\eta\trianglelefteq\nu\}.
  \end{array}\]

\item The Cantor space $\can$ and the Baire space $\baire$ are the spaces of
all functions from $\omega$ to $2$, $\omega$, respectively, equipped with
the natural (Polish) topology. 

\item The quantifiers $(\forall^\infty n)$ and $(\exists^\infty n)$ are
abbreviations for  
\[(\exists m\in\omega)(\forall n>m)\quad\mbox{ and }\quad(\forall m\in\omega)
(\exists n>m),\] 
respectively. For $f,g\in\baire$ we write $f<^* g$ ($f\leq^* g$,
respectively) whenever $(\forall^\infty n\in\omega)(f(n)<g(n))$
($(\forall^\infty n\in \omega)(f(n)\leq g(n))$, respectively). 

\item $\mbR^{{\geq}0}$ stands for the set of non-negative reals. 
\end{enumerate}

\noindent{\bf Basic convention:} In this paper, $\bH$ is a function from
$\omega$ to $\omega\setminus 2$ and $\cX=\prod\limits_{i<\omega}\bH(i)$. The
space $\cX$ is equipped with natural (Polish) product topology. 

\section{Regular universality parameters}
Since our main result applies to a somewhat restricted class of universal
parameters of \cite[\S 2.3]{RoSh:672}, we adopt here a simplified version of 
the definition of universality parameters (it fits better the case we
cover). The main difference between our definition \ref{univpar} and
\cite[Def. 2.3.3]{RoSh:672} is that we work in the setting of complete tree
creating pairs (so we may ignore $(K,\Sigma)$ and just work with trees) and
$\cF^\gp$ is assumed to be a singleton (so we also ignore it incorporating
its function into $\cG^\gp$). This simplification should increase clarity,
but we still include particular examples from \cite[\S 2.4]{RoSh:672} (see
\ref{PPex}, \ref{ucmz} at the end of this section).     

\begin{definition}
\label{Htree}
\begin{enumerate}
\item {\em A finite $\bH$--tree} is a tree $S\subseteq
\bigcup\limits_{n\leq N}\prod\limits_{i<n}\bH(i)$ with $N<\omega$,
$\mrot(S)=\langle\rangle$ and $\max(S)\subseteq\prod\limits_{i<N}\bH(i)$. 
The integer $N$ may be called {\em the level of the tree $S$} and it will be
denoted by $\lev(S)$.
\item {\em An infinite $\bH$--tree} is a tree $T\subseteq
\bigcup\limits_{n<\omega}\prod\limits_{i<n}\bH(i)$ with $\mrot(T)=\langle
\rangle$ and $\max(T)=\emptyset$.
\item The family of all finite $\bH$--trees will be denoted by $\FC[\bH]$,
  and the set of all infinite $\bH$--trees will be called $\IFC[\bH]$
\end{enumerate}
\end{definition}

\begin{definition}
[Compare {\cite[Def. 2.3.3]{RoSh:672}}] 
\label{univpar}
{\em A simplified universality parameter $\gp$ for $\bH$} is a pair
$(\cG^\gp,F^\gp)=(\cG,F)$ such that   
\begin{enumerate}
\item[$(\alpha)$] elements of $\cG$ are triples $(S,n_\dn,n_\up)$ such that 
$S$ is a finite $\bH$--tree and $n_\dn\leq n_\up\leq\lev(S)$, $(\{\langle
\rangle\},0,0)\in\cG$;  
\item[$(\beta)$] {\bf if:}\qquad $(S^0,n^0_\dn,n^0_\up)\in\cG$, $S^1$ is a
finite $\bH$--tree, $\lev(S^0)\leq\lev(S^1)$, and $S^1\cap 
\prod\limits_{i<\lev(S^0)}\bH(i)\subseteq S^0$, and $n^1_\dn\leq n^0_\dn$,
$n^0_\up\leq n^1_\up\leq\lev(S^1)$,\\  
{\bf then:}\quad $(S^1,n^1_\dn,n^1_\up)\in\cG$, 
\item[$(\gamma)$] $F\in\baire$ is increasing, 
\item[$(\delta)$] {\bf if:}
\begin{itemize}
\item $(S^\ell,n^\ell_\dn,n^\ell_\up)\in\cG$ (for $\ell<2$),
  $\lev(S^0)=\lev(S^1)$,   
\item $S\in\FC[\bH]$, $\lev(S)<\lev(S^\ell)$, and  $S^\ell\cap 
\prod\limits_{i<\lev(S)}\bH(i)\subseteq S$ (for $\ell<2$),   
\item $\lev(S)<n^0_\dn$, $n^0_\up<n^1_\dn$, $F(n^1_\up)<\lev(S^1)$,
\end{itemize}
{\bf then:}\quad there is $(S^*,n^*_\dn,n^*_\up)\in\cG$ such that
\begin{itemize} 
\item $n^*_\dn=n^0_\dn$, $n^*_\up=F(n^1_\up)$, $\lev(S^*)=\lev(S^0)=
\lev(S^1)$, and 
\item $S^0\cup S^1\subseteq S^*$ and $S^*\cap\prod\limits_{i<\lev(S)}\bH(i)=S$. 
\end{itemize}
\end{enumerate}
\end{definition}

\begin{definition}
[Compare {\cite[Def. 2.3.5]{RoSh:672}}] 
\label{unforc}
Let $\gp=(\cG,F)$ be a simplified universality parameter for $\bH$. 
\begin{enumerate}
\item We say that an infinite $\bH$--tree $T$ is {\em $\gp$--narrow} if 
for infinitely many $n<\omega$, for some $n=n_\dn<n_\up$ we have
\[(T\cap\bigcup\limits_{n\leq n_\up+1}\prod\limits_{i<n}\bH(i),
n_\dn,n_\up)\in\cG.\]
\item We define a forcing notion $\bQp$:

\noindent {\bf A condition} in $\bQp$ is a pair $p=(N^p,T^p)$ such that
$N^p<\omega$ and $T^p$ is an infinite $\gp$--narrow $\bH$--tree. 

\noindent{\bf The order $\leq$} on $\bQp$ is given by:\\
$(N^0,T^0)\leq (N^1,T^1)$\quad if and only if
\begin{itemize}
\item $N^0\leq N^1$, $T^0\subseteq T^1$, and
\item $T^1\cap\prod\limits_{i<N^0}\bH(i)= T^0\cap\prod\limits_{i<N^0}\bH(i)$. 
\end{itemize}
\end{enumerate}
\end{definition}

\begin{proposition}
[Compare {\cite[Prop. 2.3.6]{RoSh:672}}] 
\label{unccc}
If $\gp$ is a simplified universality parameter, then $\bQp$ is a Borel 
$\sigma$--centered forcing notion.
\end{proposition}

\begin{definition}
[Compare {\cite[Def. 3.2.1]{RoSh:672}}] 
\label{unideal}
Let $\gp=(\cG,F)$ be a simplified universality parameter for $\bH$.
\begin{enumerate}
\item We say that {\em $\gp$ is suitable} whenever: 
\begin{enumerate}
\item[(a)] for every $n<\omega$, there is $N>n$ such that 

\noindent {\bf if} $(S,n_\dn,n_\up)\in\cG$, $N\leq n_\dn$ and $\eta\in
\prod\limits_{i<n}\bH(i)$,  

\noindent {\bf then} $(\exists\nu\in\prod\limits_{i<\lev(S)}\bH(i))(\eta
\vartriangleleft\nu\ \&\ \nu\notin S)$, and  
\item[(b)] for every $n<\omega$, there is $N>n$ such that 

\noindent {\bf if} $S$ is a finite $\bH$--tree, $\lev(S)=n$, $\eta\in
\prod\limits_{i<N}\bH(i)$ and $\eta\rest n\in S$, 

\noindent {\bf then} there is $(S^*,n_\dn,n_\up)\in\cG$ such that
$n<n_\dn\leq n_\up<N$, $S\subseteq S^*$, $S^*\cap\prod\limits_{i<\lev(S)}
\bH(i)=\max(S)$ and $\eta\in S^*$.
\end{enumerate}
\item We say that a closed set $A\subseteq\cX$ is $\gp$--narrow if the
corresponding infinite $\bH$--tree $T$ (i.e.,  $A=[T]$) is $\gp$--narrow. 
\item $\cI_\gp^0$ is the ideal generated by $\gp$--narrow closed subsets of
  $\cX$. 
\item $\cI_\gp$ is the $\sigma$--ideal generated by $\cI^0_\gp$. 
\item $\name{T}_\gp$ is a $\bQp$--name such that 
\[\forces_{\bQp}\name{T}_\gp=\bigcup\{T^p\cap\prod_{i<N^p}\bH(i): p\in
\name{G}_\bQp\}.\]
\end{enumerate}
\end{definition}

\begin{proposition}
[Compare {\cite[Prop. 3.2.3]{RoSh:672}}]
\label{proideal} 
Let $\gp$ be a suitable simplified universality parameter for $\bH$. 
\begin{enumerate}
\item Every set in $\cI^0_\gp$ is nowhere dense in $\cX$; all singletons
belong to $\cI^0_\gp$. 
\item If $T_0,T_1\in\IFC[\bH]$ are $\gp$--narrow, then $T_0\cup T_1\in
\IFC[\bH]$ is $\gp$--narrow.  
\item $\cI^0_\gp$ is an ideal and $\cI_\gp$ is a proper Borel
$\sigma$--ideal of subsets of $\cX$. 
\item In $\bV^\bQp$, $\name{T}_\gp$ is an infinite $\gp$--narrow
  $\bH$--tree.  
\end{enumerate}
\end{proposition}

Let us recall some of the examples of universality parameters from
\cite{RoSh:672}. We represent them in a somewhat modified form to fit the
simplified setting here. 

\begin{definition}
[Compare {\cite[Ex. 2.4.9]{RoSh:672}}]
\label{PPex}
Let $g\in\baire$ and $\bF:\FC[\bH]\longrightarrow\mbR^{\geq 0}$ and
$A\in\iso$.  We define $\cG_\bF^{g,A}$ as the family consisting of
$(\{\langle\rangle\},0,0)$ and of all triples $(S,n_\dn,n_\up)$ such that  
\begin{enumerate}
\item[$(\alpha)$] $S$ is a finite $\bH$--tree, $n_\dn\leq n_\up\leq\lev(S)$,
  and 
\item[$(\beta)$] $\big(\forall\nu\in S\cap\prod\limits_{i<n_{dn}}\bH(i)
\big)\big(\exists\eta\in\prod\limits_{i<\lev(S)}\bH(i)\big)\big(\nu
\vartriangleleft\eta\ \&\ \eta\notin S\big)$,
\end{enumerate}
and such that for some sequence $\langle Y_i:i\in A\cap[n_\dn, n_\up]
\rangle$ we have  
\begin{enumerate}
\item[$(\gamma)$] $Y_i\in\FC[\bH]$, $\lev(Y_i)=i+1$, $\bF(Y_i)\leq g(i)$
(for all $i\in A\cap[n_\dn, n_\up)$), and  
\item[$(\delta)$] $\big(\forall\eta\in \max(S)\big)\big(\exists i\in A\cap 
[n_\dn,n_\up)\big)\big(\eta\rest i\in\max(Y_i)\big)$.  
\end{enumerate}
If $A=\omega$ then we may omit it and write $\cG_\bF^g$.
\end{definition}

\begin{proposition}
\label{PPprop}
Let $F_\bH\in\baire$ be an increasing function such that 
\[\big(\forall n<\omega\big)\big((n+1)^2\cdot\prod\limits_{i\leq
  n}\bH(i)<F_\bH(n)\big)\]
and let $g\in\baire$, $A\in\iso$. If a function $\bF:\FC[\bH]\longrightarrow 
\mbR^{{\geq}0}$ satisfies 
\[\big(\forall S\in\FC[\bH]\big)\big(|\max(S)|=1\ \Rightarrow\
\bF(S)=0\big),\] 
then $(\cG^{g,A}_\bF,F_\bH)$ is a suitable simplified universality parameter. 
\end{proposition}

\begin{example}
\label{PPpart}
Let $g\in\baire$, $A\in\iso$.
\begin{enumerate}
\item Let $\bF_0,\bF_1:\FC[\bH]\longrightarrow\mbR^{\geq 0}$ be defined by
\[\bF_0(S)=\max\big(|\suc_S(s)|:s\in S\setminus\max(S)\big)-1\quad\mbox{ and
}\quad\bF_1(S)=|\max(S)|-1\]
(for $S\in\FC[\bH]$). Then both $(\cG^{g,A}_{\bF_0},F_\bH)$ and
$(\cG^{g,A}_{\bF_1},F_\bH)$ are suitable simplified universality parameters.  
\item Let $\bF_2:\FC[\bH]\longrightarrow\mbR^{\geq 0}$ be defined by $\bF_2(
\{\langle\rangle\})=0$ and 
\[\bF_2(S)=\Big|\big\{\eta\big(\lev(S)-1\big):\eta\in\max(S)\big\}\Big|-1\]
when $\lev(S)>0$. Then $(\cG^{g,A}_{\bF_2},F_\bH)$ is s suitable simplified
universality parameter.  
\item Suppose that $(K,\Sigma)$ is a local tree creating pair for $\bH$ (see
  \cite[\S 1.3, Def. 1.4.3]{RoSh:470}) such that 
\begin{itemize}
\item for each $n<\omega$, $\eta\in\prod\limits_{i<n}\bH(i)$ and a non-empty
set $X\subseteq \bH(n)$, there is a unique tree creature $t_{\eta,X}\in K$
satisfying $\pos(t_{\eta,X})=\{\eta\conc\langle k\rangle:k\in X\}$,   
\item if $n<\omega$, $\eta\in\prod\limits_{i<n}\bH(i)$, $X\subseteq \bH(n)$
  and $|X|=1$, then $\nor[t_{\eta,X}]=0$.  
\end{itemize}
For $S\in\FC[\bH]$ let 
\[\bF_3(S)=\bF^{K,\Sigma}_3(S)\stackrel{\rm def}{=}\max(\nor[t_\eta]:\eta\in
\hat{S}),\]
where $\langle t_\eta:\eta\in\hat{S}\rangle$ is the unique finite
tree--candidate such that $\pos(t_\eta)=\suc_S(\eta)$ for $\eta\in\hat{S}=
S\setminus\max(S)$. Then $(\cG^{g,A}_{\bF_3},F_\bH)$ is a suitable simplified
universality parameter. 
\end{enumerate}
\end{example}

\begin{remark}
\label{easrem}
The universality parameters from \ref{PPpart} are related to the
PP--property and the strong PP--property (see \cite[Ch VI, 2.12*]{Sh:f},  
compare also with \cite[\S 7.2]{RoSh:470}). Note that if $A\in\iso$ and
$g\in\baire$ then an infinite $\bH$--tree $T$ is $(\cG^{g,A}_{\bF_2},
F_\bH)$--narrow if and only if there exist sequences $\bar{w}=\langle
w_i:i\in A\rangle$ and $\bar{n}=\langle n_k:k<\omega\rangle$ such that  
\begin{itemize}
\item $\big(\forall i\in A\big)\big(w_i\subseteq\bH(i)\ \&\ |w_i|\leq g(i)+1 
\big)$, and 
\item $n_k<n_{k+1}<\omega$ for each $k<\omega$, and 
\item $\big(\forall\eta\in [T]\big)\big(\forall k<\omega\big)\big(\exists
  i\in A\cap [n_k,n_{k+1})\big)\big(\eta(i)\in w_i\big)$.
\end{itemize}
\end{remark}

\begin{definition}
[Compare {\cite[Ex. 2.4]{RoSh:672}}]
\label{ucmz}
Let $\cG^{\rm cmz}_\bH$ consist of $(\{\langle\rangle\},0,0)$ and of all
triples $(S,n_\dn,n_\up)$ such that $S\in\FC[\bH]$, $n_\dn\leq n_\up\leq
\lev(S)$ and
\[\frac{|S\cap\prod\limits_{i<n_\up}\bH(i)|}{|\prod\limits_{i<n_\up}\bH(i)|}
\leq \sum_{i=n_\dn}^{n_\up}\frac{1}{(i+1)^2}.\]
\end{definition}

\begin{proposition}
\label{propumz}
\begin{enumerate}
\item Let $F_\bH$ be as in \ref{PPprop}. Then $\gp^{\rm cmz}_\bH=(\cG^{\rm
  cmz}_\bH,F_\bH)$ is a suitable simplified universality parameter.  
\item An infinite $\bH$--tree $T$ is $\gp^{\rm cmz}_\bH$--narrow\quad if and 
only if\quad $[T]$ is of measure zero (with respect to the product measure
on $\cX$). 
\item $\cI_{\gp^{\rm cmz}_\bH}$ is the $\sigma$--ideal of subsets of
$\cX$ generated by closed measure zero sets.
\end{enumerate}
\end{proposition}

\begin{definition}
\label{permute}
\begin{enumerate}
\item {\em A coordinate-wise permutation for $\bH$} is a sequence 
$\bar{\pi}=\langle\pi_n:n<\omega\rangle$ such that (for each $n<\omega$)
$\pi_n:\bH(n)\longrightarrow\bH(n)$ is a bijection. We say that such
$\bar{\pi}$ is {\em an $n$--coordinate-wise permutation} if $\pi_i$ is the
identity for all $i>n$. 
\item {\em A rational permutation for $\bH$\/} is an $n$--coordinate-wise
permutation for $\bH$ (for some $n<\omega$). The set of all
$n$--coordinate-wise permutations for $\bH$ will be called $\rpHn$ and the
set of all rational permutation will be denoted by $\rpH$ (so
$\rpH=\bigcup\limits_{n\in\omega}\rpHn$).  
\item Let $\bar{\pi}$ be a coordinate-wise permutation for $\bH$. We will
treat $\bar{\pi}$ as a bijection from $\bigcup\limits_{n\leq\omega}
\prod\limits_{i<n}\bH(i)$ onto $\bigcup\limits_{n\leq\omega}
\prod\limits_{i<n}\bH(i)$ such that for $\eta\in \prod\limits_{i<n}\bH(i)$
($n\leq\omega$) and $i<n$ we have $\bar{\pi}(\eta)(i)=\pi_i(\eta(i))$. 
\end{enumerate}
\end{definition}

\begin{definition}
\label{regular}
A simplified universality parameter $\gp=(\cG,F)$ for $\bH$ is called {\em a
  regular universality parameter} whenever
\begin{enumerate}
\item[(a)] $\gp$ is suitable (see \ref{unideal}(1)), and 
\item[(b)] $\cG$ is invariant under rational permutations, that is 
\begin{quote}
if $\bar{\pi}\in\rpH$ and $(S,n_\dn,n_\up)\in\cG$, then $(\bar{\pi}[S],n_\dn,
  n_\up)\in\cG$.  
\end{quote}
\end{enumerate}
\end{definition}

\begin{proposition}
\label{XPP}
\begin{enumerate}
\item Suppose that $F_\bH\in\baire$, $g\in\baire$, $A\in\iso$ and a function
  $\bF:\FC[\bH]\longrightarrow\mbR^{{\geq}0}$ are as in \ref{PPprop}. Assume
  also that 
\[\big(\forall S\in\FC[\bH]\big)\big(\forall\bar{\pi}\in\rpH\big)
  \big(\bF(S)=\bF(\bar{\pi}[S])\big).\]
Then $(\cG^{g,A}_\bF,F_\bH)$ is a regular universality parameter. 
\item For $i=0,1,2$, $(\cG^{g,A}_{\bF_i},F_\bH)$ (defined in
  \ref{PPex}(1,2)) is a regular universality parameter. 
\end{enumerate}
\end{proposition}

From now on we will assume that all universality parameters we consider are
regular. The ideals associated with regular parameters are much nicer than
those in the general case, and they are more directly connected with the
respective universal forcing notions.  

\begin{lemma}
\label{lemeasy}
Suppose that $T\in\IFC[\bH]$ is a $\gp$--narrow tree. Then there are a
$\gp$--narrow tree $T^*\in\IFC[\bH]$ and a strictly increasing sequence
$\bar{n}=\langle n_k:k<\omega\rangle\subseteq\omega$ such that 
\begin{enumerate}
\item[(a)] $T\subseteq T^*$ and for every $k<\omega$:
\item[(b)$_k$] if $\nu_0,\nu_1\in T^*\cap\prod\limits_{i\leq n_k}\bH(i)$ and 
$\bar{\pi}\in{\rm rp}^{n_k}_{\bH}$ are such that $\bar{\pi}(\nu_0)=\nu_1$,
then $\bar{\pi}[(T^*)^{[\nu_0]}]=(T^*)^{[\nu_1]}$, and   
\item[(c)$_k$] if a finite $\bH$--tree $S\in\FC[\bH]$ is such that  
\begin{itemize}
\item $\lev(S)=n_{k+1}+1$, and 
\item for all $\nu_0\in S\cap\prod\limits_{i\leq n_k}\bH(i)$ and $\nu_1\in
T^*\cap \prod\limits_{i\leq n_k}\bH(i)$ and $\bar{\pi}\in {\rm
  rp}^{n_k}_{\bH(i)}$ such that $\bar{\pi}(\nu_0)=\nu_1$ we have:
$\bar{\pi}[S^{[\nu_0]}]\subseteq (T^*)^{[\nu_1]}$,
\end{itemize}
then $(S, n_k+1, n_{k+1})\in\cG$.
\end{enumerate}
\end{lemma}

\begin{proof}
We will define $n_k$ and $T^*\cap\bigcup\limits_{n\leq n_k+1}
\prod\limits_{i<n}\bH(i)$  inductively. We let $n_0=0$ and $T^*\cap
\prod\limits_{i\leq 0}\bH(i)=T\cap \prod\limits_{i\leq 0}\bH(i)$. Now
suppose that $n_k$ and $T^*\cap\bigcup\limits_{n\leq n_k+1}
\prod\limits_{i<n}\bH(i)$ have been already chosen. Let  
\[T^+=\bigcup\{\bar{\pi}[T]:\bar{\pi}\in{\rm rp}^{n_k}_{\bH}\}.\]
It follows from \ref{proideal}(2) that $T^+$ is a $\gp$--narrow tree, so
we may pick $n_{k+1}>n_k$ such that 
\[\Big(T^+\cap\bigcup\limits_{n\leq n_{k+1}}\prod\limits_{i<n}\bH(i),n_k+1,
n_{k+1}\Big)\in\cG.\]
We choose $T^*\cap\bigcup\limits_{n\leq n_{k+1}+1}\prod\limits_{i<n} \bH(i)$
so that 
\[T^*\cap\prod\limits_{i\leq n_{k+1}}\bH(i)=\Big\{\eta\in T^+:\lh(\eta)=
n_{k+1}+1\ \&\ \eta\rest (n_k+1)\in T^*\cap\prod\limits_{i\leq n_k}\bH(i)
\Big\},\] 
completing the inductive definition. Now it should be clear that $\bar{n}$
and $T^*$ are as required.
\end{proof}

\begin{proposition}
\label{regbas}
\begin{enumerate}
\item The ideal $\cI^0_\gp$ is invariant under coordinate-wise permutations.
\item For every $A\in\cI_\gp$ there is $A^*\in\cI_\gp^0$ such that 
\[A\subseteq \bigcup\{\bar{\pi}[A^*]:\bar{\pi}\in\rpH\}.\]
\item $\forces_{\bQp}$`` $\name{T}_\gp$ is a $\gp$--narrow tree such that
for every closed set $A\in\cI_\gp^0$ coded in $\bV$, there is $n<\omega$
with $A\subseteq\bigcup\{[\bar{\pi}[\name{T}_\gp]]:\bar{\pi}\in\rpHn\}$''. 
\end{enumerate}
\end{proposition}

\begin{proof}
(3)\quad It follows from \ref{proideal}(4) that $\forces_{\bQp}$``
$\name{T}_\gp$ is $\gp$--narrow ''. 

Suppose now that $p=(N,T)\in\bQp$ and $A\subseteq\cX$ is a closed set from
$\cI^0_\gp$. Pick $S\in\IFC[\bH]$ such that  
\begin{itemize}
\item $|S\cap\prod\limits_{i\leq N}\bH(i)|=1$, say $S\cap\prod\limits_{i\leq 
  N}\bH(i)=\{\nu_0\}$, and 
\item $S$ is $\gp$--narrow, and
\item $A\subseteq\bigcup\{[\bar{\pi}[S]]:\bar{\pi}\in{\rm rp}^N_{\bH}\}$.  
\end{itemize}
Now we may pick a condition $q\in\bQp$ stronger than $p$ and such that
$N^q=N$ and 
\[\bigcup\{\bar{\pi}[S]:\bar{\pi}\in{\rm rp}^N_{\bH}\ \&\ \bar{\pi}(\nu_0)
  \in T\}\subseteq T^q.\]
Then $q\forces_{\bQp}$`` $A\subseteq\bigcup\{[\bar{\pi}[\name{T}_\gp]]:
\bar{\pi}\in {\rm rp}^N_{\bH}\}$ ''.
\end{proof}

\section{Generic objects for regular universal forcing notions}
In this section we present our main result: a sequence 
\begin{quote}
{\em Cohen real --- dominating real --- Cohen real} 
\end{quote}
produces generic filters for forcing notions $\bQp$ determined by regular
universality parameters $\gp$.

\begin{theorem}
\label{cohdomcoh}
Let $\gp=(\cG,F)$ be a regular universality parameter for $\bH$.
\begin{enumerate}
\item Suppose that $\bV\subseteq\bV^*\subseteq\bV^{**}$ are universes of set
theory, $\gp\in\bV$, $T\in\bV^*$ and $c\in\baire\cap\bV^{**}$ are such that 
\begin{enumerate}
\item[(a)] $T\in\IFC[\bH]$ is a $\gp$--narrow tree such that for every
closed set $A\in\cI_\gp^0$ coded in $\bV$, there is $n<\omega$ with
$A\subseteq \bigcup\{[\bar{\pi}[T]]:\bar{\pi}\in\rpHn\}$, and
\item[(b)] $c$ is a Cohen real over $\bV^*$.
\end{enumerate}
Then, in $\bV^{**}$, there is a generic filter $G\subseteq
\big(\bQp\big)^\bV$ over $\bV$.  
\item Suppose that $\bV\subseteq\bV^*\subseteq\bV^{**}$ are universes of set
theory, $\gp\in\bV$, $c\in\baire\cap\bV^*$ and $d\in\baire\cap\bV^{**}$
are such that  
\begin{enumerate}
\item[(a)] $c$ is a Cohen real over $\bV$, and
\item[(b)] $d$ is dominating over $\bV^*$.
\end{enumerate}
Then, in $\bV^{**}$, there is a $\gp$--narrow tree $T\in\IFC[\bH]$ such that
for every closed set $A\in\cI_\gp^0$ coded in $\bV$, there is $n<\omega$ with
$A\subseteq \bigcup\{[\bar{\pi}[T]]:\bar{\pi}\in\rpHn\}$. 
\end{enumerate}
\end{theorem}

\begin{proof}
(1)\quad The proof essentially follows the lines of that of this result
for the case of the Universal Meager forcing notion by Truss \cite[Lemma
  6.4]{Tr77}. So suppose that $T$, $c$ are as in the assumptions. Let 
$\bar{n}=\langle n_k:k<\omega\rangle,T^*\in\bV^*$ be as given by
\ref{lemeasy} for $T$ (so they satisfy \ref{lemeasy}(a--c)). 

Consider the following forcing notion $\bbC^*=\bbC^*(\bar{n},T^*)$:

\noindent {\bf A condition} in $\bbC^*$ is a finite $\bH$--tree $S$ such that 
$\lev(S)=n_k+1$ for some $k<\omega$.

\noindent{\bf The order relation $\leq_{\bbC^*}$} on $\bbC^*$ is given by:\\ 
$S_0\leq_{\bbC^*} S_1$ if and only if $S_0\subseteq S_1$ and $S_1\cap
\prod\limits_{i<\lev(S_0)}\bH(i)=\max(S_0)$, and    
\begin{enumerate}
\item[$(\otimes)$] {\em if\/} $\lev(S_0)=n_k+1$, $\lev(S_1)=n_\ell+1$,
$\nu_0\in \max(S_0)$, $\nu_1\in T^*\cap\prod\limits_{i\leq n_k}\bH(i)$ and
$\bar{\pi}\in {\rm rp}^{n_k}_{\bH}$ are such that
  $\bar{\pi}(\nu_1)=\nu_0$,\\ 
{\em then\/} $\bar{\pi}[(T^*)^{[\nu_1]}\cap\prod\limits_{i\leq
    n_\ell}\bH(i)]\subseteq (S_1)^{[\nu_0]}$.  
\end{enumerate}
Plainly, $\bbC^*$ is a countable atomless forcing notion, so it is
equivalent to the Cohen forcing $\bbC$. Therefore the Cohen real $c$
determines a generic filter $G^c\subseteq\bbC^*$ over $\bV^*$. Letting 
$T^c=\bigcup G^c$ we get an infinite $\bH$--tree, $T^c\in\bV^{**}$. By an
easy density argument, for infinitely many $k<\omega$, for each $\nu_0\in
T^c\cap\prod\limits_{i\leq n_k}\bH(i)$, $\nu_1\in T^*\cap
\prod\limits_{i\leq n_k}\bH(i)$ and $\bar{\pi}\in {\rm rp}^{n_k}_{\bH}$ such
that $\bar{\pi}(\nu_0)=\nu_1$ we have  
\[(T^*)^{[\nu_1]}\cap \prod\limits_{i\leq n_{k+1}}\bH(i)=
\pi[(T^c)^{[\nu_0]}]\cap \prod\limits_{i\leq n_{k+1}}\bH(i).\] 
Hence $T^c$ is $\gp$--narrow (remember \ref{lemeasy}(c)). Also, because of
the definition of the order,  
\begin{enumerate}
\item[$(\circledast)$] if $S\in G^c$, $\lev(S)=n_k+1$, then for every
  $\nu_0\in\max(S)$, $\nu_1\in T^*\cap\prod\limits_{i\leq n_k}\bH(i)$ and
  $\bar{\pi}\in{\rm rp}^{n_k}_{\bH}$ such that $\bar{\pi}(\nu_1)=\nu_0$ we
  have $\bar{\pi}[(T^*)^{[\nu_1]}]\subseteq T^c$.   
\end{enumerate}

Suppose now that $D\in\bV$ is an open dense subset of $\big(\bQp\big)^\bV$. 
In $\bV^*$ we define $C_D\subseteq \bbC^*$ as the collection of all
$S\in\bbC^*$ such that for some $k<\omega$ and $T'\in\IFC[\bH]\cap\bV$ we
have     
\begin{itemize}
\item $(n_k,T')\in D$, and $T'\cap\prod\limits_{i\leq n_k}\bH(i)=\max(S)$, 
  and 
\item $T'\subseteq \bigcup\{\bar{\pi}[T^*]:\bar{\pi}\in{\rm
  pr}^{n_k}_{\bH}\}$. 
\end{itemize}

\begin{claim}
\label{cl3}
$C_D\in\bV^*$ is an open dense subset of $\bbC^*$.
\end{claim}

\begin{proof}[Proof of the Claim]
Working in $\bV^*$, let $S_0\in\bbC^*$, $n_{k_0}=\lev(S_0)-1$. Pick an
infinite $\bH$--tree $S^+\in\IFC[\bH]$ such that $S^+\cap\prod\limits_{i
\leq n_{k_0}}\bH(i)=\max(S_0)$ and 
\begin{enumerate}
\item[{\bf if}] $\eta\in S^+$, $\nu\in T^*$ and $\lh(\eta)=\lh(\nu)
=n_{k_0}+1$ and $\bar{\pi}\in{\rm rp}^{n_{k_0}}_{\bH}$ are such that
$\bar{\pi}(\nu)=\eta$,
\item[{\bf then}] $(S^+)^{[\eta]}=\bar{\pi}[(T^*)^{[\nu]}]$.
\end{enumerate}
In $\bV$, take a maximal antichain $\cA\subseteq D$ of $\bQp^{\bV}$ such
that $N^p>n_{k_0}$ for each $p\in\cA$. It follows from \ref{unccc} that then
also $\cA$ is a maximal antichain of $\bQp^{\bV^*}$ in $\bV^*$. Therefore
some condition $p=(N^p,T^p)\in\cA$ is compatible with $(n_{k_0}+1,S^+)\in
\big(\bQp\big)^{\bV^*}$. Note that then ($N^p>n_{k_0}$ and) 
\[T^p\cap\prod\limits_{i\leq n_{k_0}}\bH(i)=S^+\cap\prod\limits_{i\leq
n_{k_0}}\bH(i)= \max(S_0)\]
and $S^+\cap \prod\limits_{i<N^p}\bH(i)\subseteq T^p\cap\prod\limits_{i< 
N^p}\bH(i)$. Take $k<\omega$ such that $n_k>N^p>n_{k_0}$ and $T^p\subseteq
\bigcup\{\bar{\pi}[T]:\bar{\pi}\in{\rm rp}^{n_k}_{\bH}\}$ (remember the
assumption \ref{cohdomcoh}(1)(a) on $T$). Let 
\[S_1\stackrel{\rm def}{=}(S^+\cup T^p)\cap\bigcup\limits_{n\leq n_k+1}
\prod\limits_{i<n}\bH(i)\in \FC[\bH].\]
Then $S_1\in\bbC^*$ is a condition stronger than $S_0$. 

Note that $T^p,S_1\in\bV$, so we may find $T'\in\IFC[\bH]\cap\bV$ such that: 
\begin{itemize}
\item $T'\cap\prod\limits_{i\leq n_k}\bH(i)=\max(S_1)$, and 
\item if $\eta\in T^p\cap\prod\limits_{i\leq n_k}\bH(i)$, then
$(T')^{[\eta]}= (T^p)^{[\eta]}$, and 
\item if $\eta\in (T'\setminus T^p)\cap\prod\limits_{i\leq n_k}\bH(i)$, then  
$(T')^{[\eta]}=\bar{\pi}[(T^p)^{[\nu]}]$ for some $\nu\in T^p\cap
\prod\limits_{i\leq n_k}\bH(i)$ and $\bar{\pi}\in{\rm rp}^{n_k}_{\bH}$ such
that $\bar{\pi}(\nu)=\eta$. 
\end{itemize}
It follows from the choice of $k$ and from the choice of $T^*$ (remember 
\ref{lemeasy}(a,b)) that for each $\nu_0\in T^p$, $\nu_1\in T^*$,
$\bar{\pi}\in{\rm rp}^{n_k}_{\bH}$ such that $\lh(\nu_0)=\lh(\nu_1)=n_k+1$
and $\bar{\pi}(\nu_1)=\nu_0$ we have $(T^p)^{[\nu_0]}\subseteq\bar{\pi}[
(T^*)^{[\nu_1]}]$. Therefore,  
\[T'\subseteq \bigcup\{\bar{\pi}[T^*]:\bar{\pi}\in{\rm pr}^{n_k}_{\bH}\}.\] 
It should also be clear that $(n_k,T')\in \big(\bQp\big)^\bV$ is stronger
than $p\in D$, and therefore it also belongs to $D$. Consequently,
$(n_k,T')$ witnesses that $S_1\in C_D$, proving the density of $C_D$.  

To show that $C_D$ is open suppose that $S_0\in C_D$, $S_0\leq_{\bbC^*}
S_1\in\bbC^*$. Let $\lev(S_0)=n_k+1$, $\lev(S_1)=n_\ell+1$ and let
$(n_k,T')$ witness that $S_0\in C_D$. By the definition of the order of
$\bbC^*$, $\bar{\pi}[(T^*)^{[\nu_1]}\cap\prod\limits_{i\leq n_\ell}
\bH(i)]\subseteq (S_1)^{[\nu_0]}$ whenever $\nu_0\in \max(S_0)$, $\nu_1\in
T^*\cap\prod\limits_{i\leq n_k}\bH(i)$ and $\bar{\pi}\in {\rm
  rp}^{n_k}_{\bH}$ is such that $\bar{\pi} (\nu_1)=\nu_0$. Since
$T'\subseteq \bigcup\{\bar{\pi}[T^*]:\bar{\pi}\in{\rm pr}^{n_k}_{\bH}\}$ and
$T'\cap\prod\limits_{i\leq n_k}\bH(i)=\max(S_0)$, we may conclude that
$T'\cap\prod\limits_{i\leq n_\ell}\bH(i)\subseteq \max(S_1)$. Consequently
we may find $T''\in\IFC[\bH]\cap\bV$ such that 
\begin{itemize}
\item $T'\subseteq T''$ and $T''\cap\prod\limits_{i\leq n_\ell}\bH(i)=
  \max(S_1)$, and 
\item $T''\cap\prod\limits_{i\leq n_k}\bH(i)=T'\cap\prod\limits_{i\leq n_k}
\bH(i)$, and 
\item $T''\subseteq\bigcup\{\bar{\pi}[T']:\bar{\pi}\in 
{\rm rp}^{n_\ell}_{\bH}\}\subseteq\bigcup\{\bar{\pi}[T^*]:\bar{\pi}\in 
{\rm rp}^{n_\ell}_{\bH}\}$. 
\end{itemize}
Then easily $(n_\ell,T'')\in \big(\bQp\big)^\bV$ is stronger than
$(n_k,T')$, so it belongs to $D$ and thus it witnesses that $S_1\in C_D$. 
\end{proof}

\begin{claim}
\label{cl2}
Let 
\[G\stackrel{\rm def}{=}\Big\{p\in\big(\bQp\big)^\bV: T^p\subseteq T^c\ \&\
T^p\cap\prod\limits_{i\leq N^p}\bH(i)=T^c\cap\prod\limits_{i\leq N^p}
\bH(i)\Big\}\in\bV^{**}.\]
Then $G$ is a generic filer in $\big(\bQp\big)^\bV$ over $\bV$.
\end{claim}

\begin{proof}[Proof of the Claim] 
By \ref{proideal}(2), $G$ is a directed subset of $\big(\bQp\big)^\bV$. We
need that $G\cap D\neq\emptyset$ for every open dense subset $D\in\bV$ of
$\big(\bQp\big)^\bV$. So let $D\in\bV$ be an open dense subset of $\big(
\bQp\big)^\bV$ and let $C_D$ be as defined before \ref{cl3}. It follows from
\ref{cl3} that $G^c\cap C_D\neq\emptyset$, say $S\in G^c\cap C_D$. Then for
some $k<\omega$ and $T'\in\IFC[\bH]\cap\bV$ we have 
\begin{itemize}
\item $(n_k,T')\in D$, and 
\item $T'\cap\prod\limits_{i\leq n_k}\bH(i)=\max(S)$, and $T'\subseteq
  \bigcup\{\bar{\pi}[T^*]:\bar{\pi}\in{\rm pr}^{n_k}_{\bH}\}$.
\end{itemize}
Now, by $(\circledast)$, we may conclude that $T'\subseteq T^c$ getting
$(n_k,T')\in G$. 
\end{proof}
\medskip

\noindent (2)\quad Suppose that $c,d$ and $\bV\subseteq\bV^*\subseteq
\bV^{**}$ are as in the assumptions. In $\bV$, consider the following
forcing notion $\bbC^{**}$: 
\smallskip

\noindent {\bf A condition} in $\bbC^{**}$ is a pair $(\bar{n},S)$ such that
\begin{enumerate}
\item[$(\alpha)$] $\bar{n}=\langle n_i:i\leq k\rangle\subseteq\omega$ is a
  strictly increasing finite sequence (so $k<\omega$),
\item[$(\beta)$]  $S\in\FC[\bH]$ is a finite $\bH$--tree such that
  $\lev(S)=n_k+1$, and for $\ell<k$: 
\item[$(\gamma)_\ell$] if $\nu_0,\nu_1\in S$, $\lh(\nu_0)=\lh(\nu_1)=n_\ell
  +1$, and $\bar{\pi}\in{\rm rp}^{n_\ell}_{\bH}$ is such that $\bar{\pi}(
  \nu_0)=\nu_1$, then $\bar{\pi}[S^{[\nu_0]}]=S^{[\nu_1]}$, and 
\item[$(\delta)_\ell$] {\bf if}
\begin{itemize}
\item $T\in\FC[\bH]$, $\lev(T)=n_{\ell+1}+1$ and 
\item for each $\nu_0\in S$, $\nu_1\in T$, $\lh(\nu_0)=\lh(\nu_1)=n_\ell+1$
  and $\bar{\pi}\in{\rm rp}^{n_\ell}_{\bH}$ such that $\bar{\pi}(\nu_0)=
  \nu_1$ we have $T^{[\nu_1]}\subseteq\bar{\pi}[S^{[\nu_0]}]$, 
\end{itemize}
{\bf then} there is $n<\omega$ such that $n_\ell+1<n\leq
F(n)<n_{\ell+1}$ and $(T,n_\ell+1,n)\in\cG$. 
\end{enumerate}

\noindent{\bf The order relation $\leq_{\bbC^{**}}$} on $\bbC^{**}$ is
essentially that of the end-extension:\\  
$(\bar{n}^0,S^0)\leq_{\bbC^{**}} (\bar{n}^1,S^1)$ if and only if
$\bar{n}^0\trianglelefteq\bar{n}^1$, $S_0\subseteq S_1$ and $S_1\cap 
\prod\limits_{i<\lev(S_0)}\bH(i)=\max(S_0)$. 
\smallskip

Since $\bbC^{**}$ is a countable atomless forcing notion, the Cohen real
$c\in\bV^*$ determines a generic filter $G^c\subseteq\bbC^{**}$ over
$\bV$, $G^c\in\bV^*$. Put  
\[\bar{n}^c=\bigcup\{\bar{n}\!:(\exists S)((\bar{n},S)\in G^c)\}\in\bV^*
\quad\mbox{and}\quad T^c=\bigcup\{S\!:(\exists\bar{n})((\bar{n},S)\in G^c)
\}\in\bV^*.\] 
Then $\bar{n}^c=\langle n_i^c:i<\omega\rangle\subseteq\omega$ is strictly
increasing and $T^c\in\IFC[\bH]$, and 
\[\Big(\forall k<\omega\Big)\Big(\big(\{\eta\in T^c:\lh(\eta)\leq n_{k+1}
+1\},n_k+1,n_{k+1}\big)\in \cG\Big)\]
Note that if $\nu_0,\nu_1\in T^c$, $\lh(\nu_0)=\lh(\nu_1)=n_k+1$ and
$\bar{\pi}\in {\rm rp}^{n_k}_{\bH}$ is such that $\bar{\pi}(\nu_0)=\nu_1$,
then $\bar{\pi}[(T^c)^{[\nu_0]}]=(T^c)^{[\nu_1]}$.  

Since, in $\bV^{**}$, there is a dominating real over $\bV^*$, we may find
$K^*=\{k^*_i:i<\omega\}\in [\omega]^{\textstyle\omega}\cap\bV^{**}$ (the
enumeration is increasing) such that 
\[(\forall K\in [\omega]^{\textstyle\omega}\cap\bV^*)(\forall^\infty i)(
|K\cap [k^*_i,k^*_{i+1})|>2).\]
Let $A$ be the set of all $\eta\in\cX=\prod\limits_{i<\omega}\bH(i)$ such
that   
\[\big(\forall i<\omega\big)\big(\exists \ell\in [k^*_i,k^*_{i+1})\big)
\big(\eta\rest (n_{\ell+1}+1)\in \bigcup\{\bar{\pi}[T^c]:\bar{\pi}\in{\rm
rp}^{n_\ell}_{\bH}\}\big).\]
Clearly, $A$ is a closed subset of $\cX$ (coded in $\bV^{**}$). Let
$T^*\in\IFC[\bH]\cap\bV^{**}$ be an infinite $\bH$--tree such that
$[T^*]=A$.    

\begin{claim}
\label{cl4}
The tree $T^*$ is $\gp$--narrow.
\end{claim}

\begin{proof}[Proof of the Claim]
Let $i<\omega$. For $\ell\in [k^*_i,k^*_{i+1})$ let  
\[T^\ell\stackrel{\rm def}{=}\big\{\nu\in\bigcup\limits_{n\leq n_{k^*_{i+1}}
+1}\prod\limits_{i<n}\bH(i):\nu\rest (n_{\ell+1}+1)\in\bigcup\{\bar{\pi}
[T^c]:\bar{\pi}\in{\rm rp}^{n_\ell}_{\bH}\}\big\}\]
and then let 
\[S^i\stackrel{\rm def}{=}\bigcup\{T^\ell:k^*_i\leq \ell<k^*_{i+1}\}\]
(so $T^\ell\in\FC[\bH]$ and also $S^i\in\FC[\bH]$). Note that for each
$\ell$ as above, by $(\delta)_\ell$, there is $n^\ell$ such that
$n_\ell+1<n^\ell\leq F(n^\ell)<n_{\ell+1}$ and $(T^\ell,n_\ell+1,
n^\ell)\in\cG$. Thus we may use repeatedly \ref{univpar}$(\delta)$ to
conclude that  $(S^i,n_{k^*_i}+1,n_{k^*_{i+1}})\in\cG$. Since 
\[T^*=\bigcap\limits_{i<\omega}\{\eta\in\bigcup\limits_{n<\omega}
\prod\limits_{j<n}\bH(j): \eta\rest (n_{k^*_{i+1}}+1)\in S^i\}\]  
we may easily finish the proof of the Claim. 
\end{proof}

\begin{claim}
\label{cl5}
For every $\gp$--narrow tree $T'\in\IFC[\bH]\cap\bV$ there is $k<\omega$
such that $T'\subseteq\bigcup\{\bar{\pi}[T^*]:\bar{\pi}\in{\rm
  rp}^{n_k}_{\bH}\}$. 
\end{claim}

\begin{proof}[Proof of the Claim]
Let $T'\in\bV$ be $\gp$--narrow. 
\smallskip

Suppose that $(\bar{n}^0,S^0)\in\bbC^{**}$, $\bar{n}^0=\langle n^0_i:i\leq
k\rangle$. Since $T^+\stackrel{\rm def}{=}\bigcup\{\bar{\pi}[T']:
\bar{\pi}\in{\rm rp}^{n_k^0}_{\bH}\}$ is also a $\gp$--narrow tree, we
may find $m>n_k^0+1$ such that 
\[(\{\eta\in T^+:\lh(\eta)\leq F(m)+3\},n_k^0+1,m)\in\cG.\] 
Let $\bar{n}^1=\bar{n}^0\conc\langle F(m)+2\rangle$, and let a finite
$\bH$--tree $S^1$ be such that 
\[\max(S^1)=\{\eta\in T^+\cap\!\prod\limits_{i<F^\cG(m)+3}\!\bH(i):
\eta\rest (n_k+1)\in S^0\}.\]
It should be clear that $(\bar{n}^1,S^1)\in\bbC^{**}$ is a condition
stronger than $(\bar{n}^0,S^0)$. 
\smallskip

Using the above considerations we may employ standard density arguments to
conclude that the set 
\[\begin{array}{ll}
K_{T'}\stackrel{\rm def}{=}\Big\{\ell<\omega:&\mbox{for all }\nu_0\in T^c, 
\nu_1\in T'\mbox{ such that }\lh(\nu_0)=\lh(\nu_1)=n_\ell+1,\\
&\bigcup\{\bar{\pi}[(T')^{[\nu_1]}\cap\prod\limits_{i\leq n_{\ell+1}}
\bH(i)]:\bar{\pi}\in{\rm rp}^{n_\ell}_{\bH}\ \&\ \bar{\pi}(\nu_1)=\nu_0\}
\subseteq T^c\Big\}
\end{array}\]
is infinite (and, of course, $K_{T'}\in\bV^*$). Therefore, by the choice of
the set $K^*\in\bV^{**}$, for some $N<\omega$ we have 
$(\forall i\geq N)(|K_{T'}\cap [k^*_i,k^*_{i+1})|>2)$. Thus, for each $i\geq
  N$ we may find $\ell\in (k^*_i,k^*_{i+1})$ such that 
\begin{itemize}
\item[{\bf if}] $\nu_0\vartriangleleft\eta\in T'$, $\nu_1\in T^c$,
  $\lh(\nu_0)=\lh(\nu_1)=n_\ell+1$, $\lh(\eta)=n_{\ell+1}+1$, and
  $\bar{\pi}\in{\rm rp}^{n_\ell}_{\bH}$ is such that $\bar{\pi}(\nu_0)=
  \nu_1$, 
\item[{\bf then}] $\bar{\pi}(\eta)\in T^c$. 
\end{itemize}
Hence we may conclude that
\[\big(\forall i\geq N\big)\big(\exists\ell\in (k^*_i,k^*_{i+1})\big)\big( 
T'\cap\prod\limits_{i\leq n_{\ell+1}+1}\bH(i)\subseteq\bigcup\{\bar{\pi}
[T^c]:\bar{\pi}\in{\rm rp}^{n_\ell}_{\bH}\}\big)\]
and therefore $T'\subseteq\bigcup\{\bar{\pi}[T^*]:\bar{\pi}\in{\rm
  rp}^{n_N}_{\bH}\}$ (just look at the choice of $T^*$). 
\end{proof}
\end{proof}

\begin{corollary}
\label{maincor}
Suppose that $\gp$ is a regular universality parameter for $\bH$ and $\bbP$
is either the Hechler forcing (i.e., standard dominating real forcing) or
the Universal Meager forcing (i.e., the amoeba for category  forcing). Then
$\bQp$ can be completely embedded into ${\rm RO}(\bbP*\name{\bbP})$.  
\end{corollary}

\section{Ideals $\cI_\gp$}
Let us recall that for an ideal $\cI$ of subsets of the space $\cX$ we
define cardinal coefficients of $\cI$ as follows:\\
the {\em additivity of $\cI$} is $\add(\cI)=\min\{|\cA|:\cA\subseteq  \cI\
\&\ \bigcup\cA\notin \cI\}$;\\ 
the {\em covering of $\cI$} is $\cov(\cI)=\min\{|\cA|:\cA\subseteq\cI\ \&\
\bigcup\cA=\cX\}$;\\
the {\em cofinality of $\cI$} is $\cof(\cI)=\min\{|\cA|:\cA\subseteq\cI\
\&\ (\forall B\in\cI)(\exists A\in\cA)(B\subseteq A)\}$;\\ 
the {\em uniformity of $\cI$} is $\non(\cI)=\min\{|A|:A\subseteq\cX\
\&\ A\notin \cI\}$.\\ 
The {\em dominating\/} and {\em unbounded\/} numbers are, respectively,
\[\begin{array}{ll}
\gd=&\min\{|\cF|:\cF\subseteq\baire\ \&\ (\forall g\in\baire)(\exists
f\in\cF)(g\leq^* f)\}\\
\gb=&\min\{|\cF|:\cF\subseteq\baire\ \&\ (\forall g\in\baire)(\exists
f\in\cF)(f\not\leq^* g)\}.
  \end{array}\]
Below, $\cM$ denotes the $\sigma$--ideal of meager subsets of $\cX$ (or of
any other Polish perfect space).  

For the rest of this section let us fixed a regular universality parameter
$\gp=(\cG,F)$ for $\bH$. 

\begin{corollary}
\label{corone}
$\add(\cM)\leq \add(\cI_\gp)$. 
\end{corollary}

\begin{proof}
It should be clear how the proof of \ref{cohdomcoh}(2) should be rewritten
to provide argument for 
\[\min\big(\gb,\cov(\cM)\big)\leq \add(\cI_\gp).\]
(Alternatively, see the proof of the dual version of this inequality in
\ref{dual} below.) By well known results of Miller and Truss we have
$\min\big(\gb,\cov(\cM) \big)=\add(\cM)$ (see \cite[Corollary
  2.2.9]{BaJu95}), so the corollary follows. 
\end{proof}

\begin{corollary}
\label{dual}
$\cof(\cI_\gp)\leq\cof(\cM)$.
\end{corollary}

\begin{proof}
By a well known result of Fremlin we have 
$\cof(\cM)=\max\big(\gd,\non(\cM)\big)$ (see \cite[Theorem
  2.2.11]{BaJu95}). Thus it is enough to show that 
\[\cof(\cI_\gp)\leq\max\big(\gd,\non(\cM)\big).\]
Let $\bbC^{**}$ be the forcing notion defined at the beginning of the proof
of \ref{cohdomcoh}(2). Let 
\[\cY\stackrel{\rm def}{=}\Big\{(\bar{n},T)\in \baire\times\IFC[\bH]:
\big(\forall k<\omega\big)\big((\bar{n}\rest (k+1),\{\eta\in T:\lh(\eta)\leq
n_k+1)\in \bbC^{**}\big)\Big\}\] 
be equipped with the natural Polish topology. Let $\kappa=\max\big(\gd,
\non(\cM)\big)$ and choose sequences $\langle K_\alpha:\alpha<\kappa\rangle$
and $\langle (\bar{n}^\alpha,T^\alpha):\alpha<\kappa\rangle$ so that 
\begin{enumerate}
\item[(i)] $K_\alpha=\{k^\alpha_i:i\in\omega\}\in\iso$ (the enumeration is
  increasing), 
\item[(ii)] $(\forall K\in\iso)(\exists\alpha<\kappa)(\forall^\infty
  i\in\omega)(|K\cap (k^\alpha_i,k^\alpha_{i+1})|>2)$, 
\item[(iii)] the set $\{(\bar{n}^\alpha,T^\alpha):\alpha<\kappa\}$ is not
  meager (in $\cY$).
\end{enumerate}
For $\alpha,\beta<\kappa$ and $N<\omega$ let
\[A^N_{\alpha,\beta}\stackrel{\rm def}{=}
\big\{\eta\in\cX: \big(\forall i\geq N\big)\big(\exists \ell\in [k^\alpha_i,
k^\alpha_{i+1})\big)\big(\eta\rest (n_{\ell+1}^\beta+1)\in\bigcup\{
\bar{\pi}[T^\beta]:\bar{\pi}\in{\rm rp}^{n^\beta_\ell}_{\bH}\}\big)\big\}.\] 
Then:
\begin{enumerate}
\item[$(*)_1$] Each $A^N_{\alpha,\beta}$ is a closed $\gp$--narrow subset of
  $\cX$. 
\end{enumerate}
[Why? See the proof of \ref{cl4}.]

\begin{enumerate}
\item[$(*)_2$] For each $\gp$--narrow tree $T'\in\IFC[\bH]$, there is
  $\beta<\kappa$ such that the set 
\[\begin{array}{ll}
K^\beta_{T'}\stackrel{\rm def}{=}\Big\{\ell<\omega:&\mbox{for all }\nu_0\in
T^\beta, \nu_1\in T'\mbox{ such that }\lh(\nu_0)=\lh(\nu_1)=n^\beta_\ell+1,\\
&\bigcup\{\bar{\pi}[(T')^{[\nu_1]}\cap\prod\limits_{i\leq n^\beta_{\ell+1}}
\bH(i)]:\bar{\pi}\in{\rm rp}^{n^\beta_\ell}_{\bH}\ \&\ \bar{\pi}(\nu_1)=\nu_0\}
\subseteq T^\beta\Big\}
\end{array}\]
is infinite.
\end{enumerate}
[Why? By (iii) and an argument similar to the one in the proof of
  \ref{cl5}.]

\begin{enumerate}
\item[$(*)_3$] For each $\gp$--narrow tree $T'\in\IFC[\bH]$ there are
  $\alpha,\beta<\kappa$ and $N<\omega$ such that $[T']\subseteq
  A^N_{\alpha,\beta}$.  
\end{enumerate}
[Why? By $(*)_2$+(ii) and an argument as in the proof of \ref{cl5}.]

Consequently, $\{A^N_{\alpha,\beta}:\alpha,\beta<\kappa\ \&\ N<\omega\}$ is
a cofinal family in $\cI^0_\gp$. Hence, by \ref{regbas}(2), 
\[\Big\{\bigcup\{\bar{\pi}[A^0_{\alpha,\beta}]:\bar{\pi}\in {\rm rp}_\bH\}:
  \alpha,\beta<\kappa\Big\}\]
is a basis of $\cI_\gp$. 
\end{proof}

\begin{proposition}
\label{addgp}
$\add(\cI_\gp)\leq\gb$ and $\gd\leq\cof(\cI_\gp)$.
\end{proposition}

\begin{proof}
Recall that $\cX=\prod\limits_{n<\omega}\bH(n)$ and $\cM=\cM(\cX)$ is the 
$\sigma$--ideal of meager subsets of $\cX$. We are going define two
functions 
\[\phi^*:\cM\longrightarrow\baire\quad\mbox{ and }\quad \phi:\baire
\longrightarrow\cI^0_\gp.\] 
First, for each $n<\omega$, pick a finite $\bH$--tree $S_n\in\FC[\bH]$ such
that 
\begin{enumerate}
\item[(a)] $\lev(S_n)>n$, $S_n\cap\prod\limits_{i<n}\bH(i)=
\prod\limits_{i<n} \bH(i)$,
\item[(b)] $(S_n,n,\lev(S_n))\in\cG$,
\end{enumerate}
and let $m_n=\lev(S_n)$. Put $M_n=\max\{m_j:j\leq n\}$ (for
$n<\omega$). Now, for $f\in\baire$ let $f^*\in\baire$ be defined by
$f^*(0)=f(0)$, $f^*(n+1)=M_{f^*(n)+f(n+1)}$, and let 
\[\phi(f)=\big\{\eta\in\cX: (\forall n<\omega)(\eta\rest m_{f^*(2n)}\in
S_{f^*(2n)})\big\}.\]
Note that, for any $f\in\baire$, $f^*$ is strictly increasing and
$\phi(f)\in\cI^0_\gp$. 

Now suppose that $B\subseteq\cX$ is meager and let $T_n\in\IFC[\bH]$ be such
that $T_n\subseteq T_{n+1}$(for $n<\omega$) and each $[T_n]$ is nowhere
dense (in $\cX$) and $B\subseteq\bigcup\limits_{n<\omega} [T_n]$. Let
$\phi^*(B)\in\baire$ be defined by letting $\phi^*(B)(0)=0$ and 
\[\begin{array}{ll}
\phi^*(B)(n{+}1)=\min\big\{k<\omega:& k>M_{\phi^*(B)(n)}\ \&\\
&(\forall\eta\in\!\!\!\!\prod\limits_{i\leq M_{\phi^*(B)(n)}}\!\!\!\!\bH(i))
(\exists\nu\in\prod\limits_{i<k}\bH(i))(\eta\vartriangleleft\notin T_{n+1})
\big\}. \end{array}\] 

\begin{claim}
\label{cl6}
If $f\in\baire$ and $B\subseteq\cX$ is meager, and $(\exists^\infty
n<\omega)(\phi^*(B)(n)<f(n))$, then $\phi(f)\setminus B\neq\emptyset$. 
\end{claim}

\begin{proof}[Proof of the Claim]
Assume that $(\exists^\infty n<\omega)(\phi^*(B)(n)<f(n))$. Then the set 
\[K=\big\{n<\omega: (\exists k<\omega)(f^*(2n)\leq\phi^*(B)(k)<\phi^*(B)
(k+1)<f^*(2n+2))\big\}\] 
is infinite. Now we may pick $\eta\in\cX$ such that for each $n<\omega$ we
have:
\begin{enumerate}
\item[(i)]  $\eta\rest m_{f^*(2n)}\in S_{f^*(2n)}$, and
\item[(ii)] if $n\in K$, then for some $k$ such that 
\[f^*(2n)\leq\phi^*(B)(k)<\phi^*(B)(k+1)<f^*(2n+2)\]
we have $\eta\rest\phi^*(B)(k+1)\notin T_{n+1}$. 
\end{enumerate}
It should be clear that the choice is possible; note that for $n,k$ as in
(ii) we have 
\[f^*(2n)<\lev(S_{f^*(2n)})=m_{f^*(2n)}<M_{f^*(2n)}\leq M_{\phi^*(B)(k)}.\] 
\end{proof}

The proposition follows from \ref{cl6}: if $\cF\subseteq\baire$ is an
unbounded family, then $\bigcup\{\phi(f):f\in\cF\}\notin\cI_\gp$, and if
$\cB\subseteq\cI_\gp$ is a basis of $\cI_\gp$, then $\{\phi^*(B):B\in\cB\}$
is a dominating family in $\baire$. 
\end{proof}

It was shown in \cite{BrSh:439} that the additivity of the $\sigma$--ideal
generated by closed measure zero sets (i.e., the one corresponding to
$\gp^{\rm cmz}_\bH$ of \ref{ucmz}) is $\add(\cM)$. We have a similar result
for another specific case of $\cI_\gp$:

\begin{proposition}
Suppose that $\bH:\omega\longrightarrow\omega\setminus 2$ and $g:\omega
\longrightarrow\omega\setminus 2$ is such that $g(n)+1<\bH(n)$ for all
$n<\omega$. Let $A\in\iso$ and $\gp=(\cG^{g,A}_{\bF_2},F_{\bH})$ (see
\ref{PPpart}(2)). Then $\add(\cI_\gp)=\add(\cM)$. 
\end{proposition}

\begin{proof}
Since $\gp$ is a regular universality parameter (by \ref{XPP}), we know that
$\add(\cM)\leq\add(\cI_\gp)\leq\gb$ (by \ref{corone}, \ref{addgp}). So for
our assertion it is enough to show that $\add(\cI_\gp)\leq\cov(\cM)$. 

Let us start with analyzing sets in $\cI_\gp$. Suppose that
$\bar{n},\bar{w}$ are such that 
\begin{enumerate}
\item[$(\otimes)_0$] $\bar{n}=\langle n_k:k<\omega\rangle$ is a strictly
  increasing sequence of integers such that $A\cap [n_k,n_{k+1})\neq
  \emptyset$ for each $k<\omega$, 
\item[$(\otimes)_1$] $\bar{w}=\langle w_i:i\in A\rangle$, $w_i\in
  [\bH(i)]^{\textstyle g(i)+1}$ for each $i\in A$.
\end{enumerate}
Put 
\[Z(\bar{n},\bar{w})\stackrel{\rm def}{=}\big\{\eta\in\prod_{i<\omega}
\bH(i):(\forall^\infty k<\omega)(\exists i\in A\cap [n_k,n_{k+1}))(\eta(i)
\in w_i)\big\}.\]
It follows from \ref{easrem} that $Z(\bar{n},\bar{w})\in\cI_\gp$. Moreover,
for every $Z\in \cI_\gp$ there are $\bar{n},\bar{w}$ satisfying $(\otimes)_0 
+(\otimes)_1$ and such that $Z\subseteq Z(\bar{n},\bar{w})$ (by
\ref{easrem}+\ref{regbas}(2)). 

\begin{claim}
\label{cl7}
Suppose that $\bar{n}^\ell,\bar{w}^\ell$ satisfy $(\otimes)_0+(\otimes)_1$
above (for $\ell=0,1$). Assume that $Z(\bar{n}^0,\bar{w}^0)\subseteq
Z(\bar{n}^1,\bar{w}^1)$. Then $(\exists^\infty k<\omega)(\forall i\in A\cap
[n^0_k,n^0_{k+1}))(w^0_i=w^1_i)$. 
\end{claim}

\begin{proof}[Proof of the Claim]
If the assertion fails, then (as $|w^0_i|=|w^1_i|=g(i)+1<\bH(i)$) we have 
$(\forall^\infty k<\omega)(\exists i\in A\cap [n^0_k,n^0_{k+1}))(w^0_i
\setminus w^1_i\neq\emptyset)$. Consequently we may pick $\eta\in
Z(\bar{n}^0,\bar{w}^0)$ such that $(\forall i\in A)(\eta(i)\notin
w^1_i)$. Then $\eta\notin Z(\bar{n}^1,\bar{w}^1)$, contradicting $Z(
\bar{n}^0,\bar{w}^0)\subseteq Z(\bar{n}^1,\bar{w}^1)$.
\end{proof}

\begin{claim}
\label{cl8}
Suppose that $f:\omega\longrightarrow\omega\setminus 2$, $\kappa<\add(
\cI_\gp)$ and $\{f_\alpha:\alpha<\kappa\}\subseteq\prod\limits_{i<\omega}
f(i)$. Then there is a function $f^*\in\prod\limits_{i<\omega} f(i)$ such
that 
\[\big(\forall\alpha<\kappa\big)\big(\exists^\infty i<\omega\big)\big(
f^*(i)=f_\alpha(i)\big).\]
\end{claim}

\begin{proof}[Proof of the Claim]
Pick an increasing sequence $\langle a_i:i<\omega\rangle$ of members of $A$
so that $f(i)<\Big|[\bH(a_i)]^{\textstyle g(a_i)+1}\Big|$ for all
$i<\omega$). For each $i$ fix a one-to-one mapping $\psi_i:f(i)
\longrightarrow [\bH(a_i)]^{\textstyle g(a_i)+1}$. Now, for $\alpha<\kappa$,
$k<\omega$ and $j\in A$ let 
\[n^\alpha_k=a_k\qquad\mbox{ and }\qquad w^\alpha_k=\left\{
\begin{array}{ll} 
\psi_i\big(f_\alpha(i)\big)&\mbox{ if }j=a_i,\ i<\omega,\\
g(j)+1&\mbox{ if }j\notin\{a_i: i<\omega\}.
\end{array}\right.\]
Then $\bar{n}^\alpha,\bar{w}^\alpha$ satisfy $(\otimes)_0+(\otimes)_1$ above
and thus $Z(\bar{n}^\alpha,\bar{w}^\alpha)\in\cI_\gp$ (for all
$\alpha<\kappa$). Since $\kappa<\add(\cI_\gp)$ we know that 
$\bigcup\limits_{\alpha<\kappa}Z(\bar{n}^\alpha,\bar{w}^\alpha)\in\cI_\gp$
and therefore we may find $\bar{n},\bar{w}$ such that they satisfy
$(\otimes)_0+(\otimes)_1$ and $\big(\forall\alpha<\kappa\big)\big(Z(
\bar{n}^\alpha,\bar{w}^\alpha)\subseteq Z(\bar{n},\bar{w})\big)$. It follows
from \ref{cl7} that 
\[\big(\forall \alpha<\kappa\big)\big(\exists^\infty k<\omega\big)\big(
\psi_k(f_\alpha(k))=w^\alpha_{a_k}=w_{a_k}\big).\] 
Let $f^*\in\prod\limits_{i<\omega} f(i)$ be such that if $k<\omega$ and
$w_k\in \rng(\psi_k)$, then $\psi_k(f^*(k))=w_k$. It should be clear that
then $f^*$ is as required.  
\end{proof}

The proposition follows now from \ref{cl8} and the inequality
$\add(\cI_\gp)\leq\gb$. 
\end{proof}

To generalize the above result to the ideals $\cI_\gp$ (for a regular
universality parameter $\gp$) one would like to know the answer to the
following question. 

\begin{problem}
Suppose that $\gp$ is a regular universality parameter for $\bH$. Does
$\add(\cI_\gp)\leq \cov(\cM)$?
\end{problem}


\begin{thebibliography}{10}

\bibitem{BaJu95}
Tomek Bartoszy\'nski and Haim Judah.
\newblock {\em {Set Theory: On the Structure of the Real Line}}.
\newblock A K Peters, Wellesley, Massachusetts, 1995.

\bibitem{BrSh:439}
Tomek Bartoszynski and Saharon Shelah.
\newblock {Closed measure zero sets}.
\newblock {\em {Annals of Pure and Applied Logic}}, 58:93--110, 1992.
\newblock math.LO/9905123.

\bibitem{J}
Thomas Jech.
\newblock {\em {Set theory}}.
\newblock Academic Press, New York, 1978.

\bibitem{JuRo}
Haim Judah and Andrzej Ros{\l}anowski.
\newblock {On Shelah's Amalgamation}.
\newblock In {\em Set Theory of the Reals}, volume~6 of {\em Israel
  Mathematical Conference Proceedings}, pages 385--414. 1992.

\bibitem{RoSh:856}
Andrzej Roslanowski and Saharon Shelah.
\newblock {How much sweetness is there in the universe?}
\newblock {\em Preprint}.
\newblock math.LO/0406612.

\bibitem{RoSh:470}
Andrzej Roslanowski and Saharon Shelah.
\newblock {Norms on possibilities I: forcing with trees and creatures}.
\newblock {\em {Memoirs of the American Mathematical Society}}, 141(671):xii +
  167, 1999.
\newblock math.LO/9807172.

\bibitem{RoSh:672}
Andrzej Roslanowski and Saharon Shelah.
\newblock {Sweet {\&} Sour and other flavours of ccc forcing notions}.
\newblock {\em Archive for Mathematical Logic}, 43:583--663, 2004.
\newblock math.LO/9909115.

\bibitem{Sh:176}
Saharon Shelah.
\newblock {Can you take Solovay's inaccessible away?}
\newblock {\em {Israel Journal of Mathematics}}, 48:1--47, 1984.

\bibitem{Sh:f}
Saharon Shelah.
\newblock {\em {Proper and improper forcing}}.
\newblock {Perspectives in Mathematical Logic}. {Springer}, 1998.

\bibitem{St85}
Jacques Stern.
\newblock Regularity properties of definable sets of reals.
\newblock {\em Annals of Pure and Applied Logic}, 29:289--324, 1985.

\bibitem{Tr77}
John Truss.
\newblock {Sets having calibre $\aleph_1$}.
\newblock In {\em {Logic Colloquium 76}}, volume~87 of {\em Studies in Logic
  and the Foundations of Mathematics}, pages 595--612. North-Holland,
  Amsterdam, 1977.

\end{thebibliography}
\def\germ{\frak} \def\scr{\cal} \ifx\documentclass\undefinedcs
  \def\bf{\fam\bffam\tenbf}\def\rm{\fam0\tenrm}\fi 
  \def\defaultdefine#1#2{\expandafter\ifx\csname#1\endcsname\relax
  \expandafter\def\csname#1\endcsname{#2}\fi} \defaultdefine{Bbb}{\bf}
  \defaultdefine{frak}{\bf} \defaultdefine{mathfrak}{\frak}
  \defaultdefine{mathbb}{\bf} \defaultdefine{mathcal}{\cal}
  \defaultdefine{beth}{BETH}\defaultdefine{cal}{\bf} \def\bbfI{{\Bbb I}}
  \def\mbox{\hbox} \def\text{\hbox} \def\om{\omega} \def\Cal#1{{\bf #1}}
  \def\pcf{pcf} \defaultdefine{cf}{cf} \defaultdefine{reals}{{\Bbb R}}
  \defaultdefine{real}{{\Bbb R}} \def\restriction{{|}} \def\club{CLUB}
  \def\w{\omega} \def\exist{\exists} \def\se{{\germ se}} \def\bb{{\bf b}}
  \def\equivalence{\equiv} \let\lt< \let\gt> \def\implies{\Rightarrow}

\end{document}